\documentclass[12pt]{article}

\usepackage{latexsym,amssymb}

\usepackage{amsmath}
\usepackage{euscript}
\usepackage{eufrak}

\newenvironment{dem}{\noindent\bf Proof. \rm}{\hfill $\mbox{\boldmath{$ \blacksquare$}}$}

\newcommand{\rr}{\vspace*{4mm}}

\newcommand{\C}{{\cal C}}
\newcommand{\LL}{{\cal L}}
\newcommand{\M}{{\cal M}}
\newcommand{\im}{\!\to\!}
\newcommand{\deb}{\!\mapsto\!}
\newcommand{\fun}{\!\longrightarrow\!}

\newtheorem{sdefi}{\bf Definition}
\newtheorem{steo}{\bf Theorem}
\newtheorem{scor}{\bf Corollary}
\newtheorem{sob}{\bf Remark}

\usepackage[english]{babel}
\begin{document}

\author{}

\title{}

\date{}

\maketitle
\thispagestyle{empty}

\begin{center}

  Marina B. LATTANZI and Alejandro G. PETROVICH

\vspace*{2.5cm}

{\large \bf A DUALITY FOR (n+1)-VALUED MV-ALGEBRAS}
\end{center}

\vspace*{2.5cm}
{\small \texttt{Abstract.} 
$MV-$algebras were introduced by Chang to prove the completeness of the infinite-valued {\L}ukasiewicz propositional calculus. 
In this paper we give a categorical equi\-va\-lence between the varieties of $(n+1)-$valued MV-algebras and the classes of Boolean algebras endowed with a certain family of filters.  Another similar categorical equi\-va\-lence is given by A. Di Nola and A. Lettieri.
 Also, we point out the relations between this categorical equivalence and the duality established by R. Cignoli, which can be derived from results obtained by P. Niederkorn on natural dualities for varieties of $MV-$algebras.}

\vspace*{2cm}
\centerline {\large \bf 1. Introduction and Preliminaries}
\vspace*{0.5cm}

Wajsberg algebras (see [6, 12, 8]) are an equivalent reformulation of Chang $MV-$ algebras based on implication instead of disjunction. $MV-$algebras were introduced by Chang [3, 4] to prove the completeness of the infinite valued {\L}ukasiewicz propositional calculus. The classes of $(n+1)-$valued $MV-$algebras were introduced by R. Grigolia in [9], who also gave their equational characterization. For each $n > 0$, this variety is generated by the chain of length $n+1$ and the algebras belonging to this variety are the algebraic models of the $(n+1)-$valued {\L}ukasiewicz propositional calculus.

Y. Komori [12] introduced the $CN-$algebras 
as algebraic models of {\L}ukasiewicz infinite-valued propositional calculus formulated in terms of the operations implication and negation. A. J. Rodriguez [15] called Wajsberg algebras what was previously known as $CN-$algebras (see also [8]). $(n+1)-$valued Wajsberg algebras are equivalent to $(n+1)-$valued $MV-$algebras. We shall deal with Wajsberg algebras instead of $MV-$algebras.

In this paper we give a direct proof of the facts that every $(n+1)-$valued Wajsberg algebra $A$ can be obtained from a boolean algebra $B$ endowed with a family of filters of $B$. The main result of this paper is proved in Theorem \ref{Teo2}, where we give an explicit construction of the isomorphism between $A$ and the Wajsberg algebra obtained in this way. This construction arises like a way to generalize the well-known results on the Post algebras of order $n+1$, which can be represented as the algebra of all functions of $\{1,2, \ldots, n\}$ into a Boolean algebra (see [1]). We also give a categorical equivalence between the varieties of $(n+1)-$valued Wajsberg algebras and the classes of pairs $\langle B,h \rangle$ where $B$ is a Boolean algebra and $h$ is a function from the lattice of divisors of $n$ into the lattice of filters of $B$ which satisfies certain conditions. 
Another similar categorical equivalence is given by A. Di Nola and A. Lettieri in [7]. Also, we point out the relations between this categorical equivalence and the duality established by R. Cignoli in [5] which can be derived from results on natural dualities for varieties of $MV-$algebras obtained by P. Niederkorn in [14]. The mentioned dualities of Cignoli and Niederkorn were derived from a general result of K. Keimel and H. Werner [11] for semiprimal varieties.

The basic results about $MV-$algebras can be found, por instance, in [6]. For a reformulation in the context of Wajsberg algebras (or $CN-$algebras) see [15, 8, 12]. 

A Wajsberg algebra (or $W-$algebra, for short) is an algebra $A=\langle A, \im , \neg, 1 \rangle$ of type $(2,1,0)$ satisfying the following identities: $\, 1\im x=x$, $\, (x\im y)\im ((y\im z)\im (x\im z))=1$, $\, (x\im y)\im y=(y\im x)\im x$  and $\, (\neg y \im\, \neg x)\im (x\im y)=1$. The reduct $(A,\vee ,\wedge ,\neg ,0,1)$ is a Kleene algebra where $0=\, \neg 1$, $x\vee y  = (x\im y)\im y$, $x\wedge y  = \,\neg (\neg x\vee \neg y)$ and $x\leq y$ if and only if $x\im y=1$. If we set $x \oplus y =\, \neg y \im x$ and $x\odot y =  \neg (x\im\, \neg y)$ then $\langle A, \oplus , \odot, 0 \rangle$ is an $MV-$algebra. The set $B(A)=\{ x \in A: x \odot x=x\}$ is a Boolean algebra; actually it is the maximal Boolean subalgebra of $A$. For all $x \in A$ and each natural number $m$ we set $x^0=1$ and $x^m=x^{m-1} \odot x$ for $m \geq 1$.

A subset $F\subseteq  A$ is an {\em implicative filter} of $A$ if $1 \in F$ and for all $a,b \in A$, $a,a\im b\in F$ implies $b\in F$.  Implicative filters are lattice filter which are closed by the operation $\odot$.  The family of all implicative filters of $A$ is an algebraic lattice under set-inclusion, and it is isomorphic to the algebraic lattice of all congruence relations on $A$. For every implicative filter $F$ of $A$ and each $x \in A$ we represent with $[x]_F$ the set of all elements $y \in A$ such that $x$ and $y$ are $F-congruent$. An implicative filter of $A$ is {\em prime} if is a lattice prime filter of $A$. We denote by $\chi(A)$ the set of all prime implicative filters of $A$. An implicative filter $P$ of $A$ is prime if and only if $A/P$ is a chain.

The unit interval $[0,1]$ endowed with the operations $x\im y:=min\ \{1,1-x+y\}$ and $ \, \neg x:=1-x$ is a $W-$algebra. For each positive integer $n$  we denote by $L_{n+1}$ the subalgebra of $[0,1]$ whose universe is $\left\{ 0,\frac{1}{n},\frac{2}{n},\ldots ,\frac{n-1}{n},1\right\}$. 
It is verified that $L_{t+1}$ is a subalgebra of $L_{n+1}$ if and only if $t$ divides $n$.

\rr

If $\langle A, \im , \neg , 1 \rangle$ is an $(n+1)-$valued Wajsberg algebra then  $\langle A, \vee , \wedge, \neg, \sigma_1,  \sigma_2, \ldots, \sigma_n, 0, 1 \rangle$ is an $(n+1)-$valued {\L}ukasiewicz algebra, where the operators $\sigma_i$, for $1 \leq i \leq n$, are defined from the $W-$operations (see [10]).

\rr

Let $B$ be a Boolean algebra and $n$ be an integer, $n \geq 1$. We denote by $B^{[n]}$ the set of all increasing monotone functions from $\{1, 2, \ldots , n\}$ into $B$. $B^{[n]}$ with the operations of the lattice defined pointwise, the chain of constants $0=c_0 < c_1 < \ldots < c_{n-1}<c_n=1$,  where $c_k(i) = \left\{ \begin{array}{rr}
                                                  1  & \mbox{if } i \geq n+1-k \\
                                                  0  & \mbox{if } i < n+1-k
                                                  \end{array}\right.$ for $0 \leq k \leq n$, the negation defined by $(\neg f)(i)= \neg f(n+1-i)$ for each $1 \leq i \leq n$ and the modal operators $\sigma_i (f)(j)=f(i)$ for all $1 \leq i \leq n$ and $1 \leq j \leq n$, is a Post algebra of order $n+1$ (see [1]); therefore it is an $(n+1)-$valued Wajsberg algebra (see  [16]). In Theorem \ref{Teo1} a direct proof of this result is given, showing explicitly the form of the operations. In every $(n+1)-$valued Wajsberg algebra, the prime filters occur in finite and disjoint chains, then by the Mart\'{\i}nez's Unicity Theorem [13] the implication is determined by the order.

\newpage

\centerline {\large \bf 2. The duality}
\vspace*{0.5cm}

\begin{steo}\label{Teo1}
Let $B$ be a Boolean algebra and $n \geq 1$ be an integer. Then $\langle B^{[n]}, \deb, \neg, \mathbb I \rangle$ is an $(n+1)-$valued Wajsberg algebra where $B^{[n]}=\{f: \{1, 2, \ldots , n\} \fun B: f(i) \leq f(j) \mbox{ for all $i,j $ such that } i \leq j \}$, $\mathbb I$ is the constant function equal to $1$ and, for $f,g \in B^{[n]}$ and $1 \leq k \leq n$, $(\neg f)(k)= \neg f(n+1-k)$ and $(f \deb g)(k)= \bigwedge\limits_{i=1}^{n-k+1} (f(i) \im g(i+k-1))$.
\end{steo}

\noindent {\bf Proof.} Let $f,g \in B^{[n]}$ and let $k,t$ be integers, $1 \leq k \leq t \leq n$ such that $(f \deb g)(k) \not\leq (f \deb g)(t)$. Then there is a prime filter $P$ of $B$ which verifies (1) $(f \deb g)(k) \in P$ and (2) $(f \deb g)(t) \notin P$. From (2) there is an integer $i_0$, $1 \leq i_0 \leq n-t+1$ such that $\neg f(i_0) \vee g(i_0 +t-1) \notin P$. This last condition implies (3) $f(i_0) \in P$ and (4) $g(i_0+t-1) \notin P$. It follows from (1) that $\neg f(i_0) \vee g(i_0 +k-1) \in P$ and by (3) we have $g(i_0 +k-1) \in P$, thus it results $g(i_0 +t-1) \in P$ which contradicts (4). Therefore the operation $\deb$ is well defined.

Let $f,g,h \in B^{[n]}$. The following properties hold:
\begin{itemize}
\item[{\rm (1)}] $\mathbb I \deb f =f$

$(\mathbb I \deb f) (k) = \bigwedge\limits_{i=1}^{n-k+1} (1 \im f(i+k-1))= \bigwedge\limits_{i=1}^{n-k+1} f(i+k-1)=f(k)$.

\item[{\rm (2)}] $(f \deb g) \deb ((g \deb h) \deb (f \deb h))= \mathbb I$

Let us observe that $f \leq g$ if and only if $f \deb g = \mathbb I$. In fact, if $f \leq g$ we have  $f(i) \leq g(i) \leq g(i+k-1)$ for all $i,k$, $1 \leq i,k \leq n$, then $f \deb g = \mathbb I$. Reciprocally, if $f \deb g = \mathbb I$ then $(f \deb g)(1) = 1$, i.e.  $f(i) \leq g(i)$ for each $i$, $1 \leq i \leq n$. 
So, it is enough to verify $(f \deb g) \leq  (g \deb h) \deb (f \deb h)$. Let $k$ be an integer, $1 \leq k \leq n$, such that $(f \deb g)(k)  \not\leq  ((g \deb h) \deb (f \deb h))(k)$; thus there is a prime filter $P$ of $B$ which verifies (5) $(f \deb g)(k) \in P$ and  (6) $((g \deb h) \deb (f \deb h))(k) \notin P$. From (5) it results (7) $f(j) \im g(j+k-1) \in P$ for any $j$, $1 \leq j \leq n-k+1$. Besides, from (6) there is an integer $j_0$, $1 \leq j_0 \leq n-k+1$ such that $(g \deb h)(j_0) \im (f \deb h) (j_0 +k -1) \notin P$, i.e., the statements (8) $(g \deb h)(j_0) \in P$ and (9) $(f \deb h) (j_0 +k -1) \notin P$ hold. It follows from (9) that there is an integer $i_0$, $1 \leq i_0 \leq n- j_0 -k+2$ which verifies (10) $f(i_0) \im h(i_0 +j_0 +k-2) \notin P$. Moreover, from (8) we obtain (11) $g(i_1) \im h(i_1 +j_0 -1) \in P$ where $i_1= i_0 +k -1 \leq n - j_0 +1$. Since $i_1 +j_0 -1= i_0 + j_0 +k -2$ from (10) and (11) we have (12) $\neg g(i_1) \in P$. On the other hand, $1 \leq j_0$ implies $n -j_0 +2 \leq n+1$, hence $i_0 +k \leq n -j_0 +2 \leq n+1$, i.e., $i_0 \leq n-k+1$; therefore from (7) we have $\neg f(i_0) \in P$ or $g(i_0 +k-1) \in P$ while from (10) and (12) it results  $f(i_0) \wedge \neg g(i_0 +k-1) \in P$  which is a contradiction.

\item[{\rm (3)}] $(f \deb g) \deb g= (g \deb f) \deb f$

We shall prove by induction on $n$ that for every integer $k$, $1 \leq k \leq n$, $((f \deb g) \deb g)(k) = f(k) \vee g(k)$ .

For $n=1$ is trivial. For each $n \geq 1$, $1 \leq k \leq n$ and $1 \leq h \leq n -k +1$ let  

$C(n,k,h) = g(h+k-1) \vee \bigvee\limits_{i=1}^{n-h+1} (f(i) \wedge \neg g(i+h-1))$ and

$U(n,k)=\bigwedge\limits_{h=1}^{n-k+1} C(n,k,h)$, i.e., $((f \deb g) \deb g)(k)=U(n,k)$.

Suppose $U(n,k)= f(k) \vee g(k)$ for all integer $k$, $1 \leq k \leq n$.

Let $k$ be an integer, $1 \leq k \leq n+1$, 

$U(n+1,k)= \bigwedge\limits_{h=1}^{n-k+2} \left( g(h+k-1) \vee \bigvee\limits_{i=1}^{n-h+2} (f(i) \wedge \neg g(i+h-1)) \right)$. Separating the term for $i=n-h+2$ we have $U(n+1,k)=$

$\bigwedge\limits_{h=1}^{n-k+2} \left(  (f(n-h+2) \wedge \neg g(n+1)) \vee  g(h+k-1) \vee \bigvee\limits_{i=1}^{n-h+1} (f(i) \wedge \neg g(i+h-1))\right)$ and  separating the term for  $h=n-k+2 $ it results

$U(n+1,k)= Z(n+1,k) \wedge \bigwedge\limits_{h=1}^{n-k+1} (  (f(n-h+2) \wedge \neg g(n+1)) \vee C(n,k,h))$ where 

$Z(n+1,k)= (f(k) \wedge \neg g(n+1)) \vee g(n+1) \vee \bigvee\limits_{i=1}^{k-1} (f(i) \wedge \neg g(i+n-k+1))=f(k) \vee g(n+1)$. Then, replacing $Z(n+1,k)$ and applying  the inductive hypothesis we obtain

$U(n+1,k)= ( f(k) \wedge \neg g(n+1)) \vee ((f(k) \vee g(n+1) \wedge (f(k) \vee g(k)) = f(k) \vee g(k)$.

\item[{\rm (4)}] $(\neg g \deb  \neg f ) \deb  (f \deb g)= \mathbb I$

It is enough to prove $(\neg g \deb \neg f)(k) \im (f \deb g)(k)=1$.

$(\neg g \deb \neg f)(k)= \bigwedge\limits_{i=1}^{n-k+1} (f(n+2-i-k) \im g(n+1-i))$, thus  

$(\neg g \deb  \neg f ) (k) \im  (f \deb g)(k)= $

$\bigwedge\limits_{j=1}^{n-k+1} \left( (\neg f(j) \vee g(j+k-1)) \vee  \bigvee\limits_{i=1}^{n-k+1}  (f(n+2-i-k) \wedge \neg g(n+1-i))  \right)= $

$\bigwedge\limits_{j=1}^{n-k+1} \left(  \bigvee\limits_{i=1}^{n+2-k-j} (\neg f(j) \vee g(j+k-1) \vee g(n+1-i) ) \vee \right.$

$\left.    \bigvee\limits_{i=n+3-k-j}^{n-k+1} (\neg f(j) \vee g(j+k-1) \vee f(n+2-i-k) ) \right)=1$.\hfill$\mbox{\boldmath{$\blacksquare$}}$
\end{itemize}

\rr

\begin{sob}\label{O1}{\rm Let $n \geq 1$ be an integer. We denote by $Div(n)$ the set of all positive divisors of $n$. Let $d \in Div(n)$. For each integer $j$, $1 \leq j \leq n$, there exists an only integer $q_{d,j}$, $1 \leq q_{d,j} \leq d$, such that $(q_{d,j} -1) \frac{n}{d} < j \leq q_{d,j}  \frac{n}{d}$. Indeed, $q_{d,j} $ is the first element of the set $X=\{q \in \mathbb N : 1 \leq q \leq d, \, j \leq q \frac{n}{d}\}$. That is to say that the only block corresponding to the divisor $d$ of $n$ that contains $j$ is that determined by $q_{d,j}$.

Thus, for any $d \in Div(n)$, we can \textit{think} an $n-$tuple to be composed by $d$ blocks, each one of them with $\frac{n}{d}$ elements.

For short, in what follows we shall write $\xi_{d,q}(f)$ instead of $f(q \frac{n}{d}) \im f((q-1) \frac{n}{d}+1)$, for each $f \in B^{[n]}$, $d \in Div(n)$ and any integer $1 \leq q \leq d$.
}
\end{sob}

\begin{scor}\label{C1}
Let $B$ be a Boolean algebra, $n \geq 1$ be an integer and $h$ be a function from the lattice of divisors of $n$ into the lattice of filters of $B$. Let $M(B,h)$ be the set \linebreak $\{f \in B^{[n]}: \xi_{d,q}(f) \in h(d) \mbox{, for each $d \in Div (n)$ and all $1 \leq q \leq d$} \}$. Then $\langle M(B,h), \linebreak \deb, \neg, \mathbb I \rangle$ is an $(n+1)-$valued Wajsberg subalgebra of $B^{[n]}$. Also, if $h(d)=B$ for each $d \in D= Div(n) -\{n\}$ then $M(B,h)$ is a Post algebra of order $n+1$.
\end{scor}

\noindent {\bf Proof.} Let $f,g \in M\left(B,h\right)$.

\begin{itemize}

\item[{\rm (a)}] $ \neg f \in M\left(B,h\right)$. In fact, given $d \in Div(n)$ and an integer  $q$, $1 \leq q \leq d$, we have  $(\neg f) (q \frac{n}{d}) \im (\neg f) ((q-1) \frac{n}{d}+1) = \neg  f(n + 1 - q \frac{n}{d}) \im \neg f( n -  (q-1) \frac{n}{d})= f(q^{\prime}  \frac{n}{d}) \im f( (q^{\prime}-1) \frac{n}{d}+1) \in F_d$, because $q^{\prime}= d-q+1$ verifies $1 \leq q^{\prime}  \leq d$.

\item[{\rm (b)}] $f \deb g  \in M\left(B, h\right)$. Indeed, let $d \in Div(n)$ and let $q$ be an integer, $1 \leq q \leq d$.

$\xi_{d,q}(f \deb g)= \bigvee\limits_{i=1}^{n-q \frac{n}{d}+1} (f(i) \wedge \neg g(i + q \frac{n}{d}-1)) \vee \bigwedge\limits_{r=1}^{n-(q-1) \frac{n}{d}} (\neg f(r) \vee g(r + ( q-1) \frac{n}{d}))= \bigwedge\limits_{r=1}^{n-(q-1) \frac{n}{d}} \ \bigvee\limits_{i=1}^{n-q \frac{n}{d}+1} (  (f(i) \vee \neg f(r) \vee  g(r + ( q-1) \frac{n}{d})) \wedge (\neg g(i + q \frac{n}{d}-1) \vee \neg f(r) \vee    g(r + ( q-1) \frac{n}{d}))  )$. 

If $r \leq \frac{n}{d}$ then $q_{d,r}=1$, thus $f(\frac{n}{d}) \im f(1) \in F_d$. Therefore  $f(\frac{n}{d}) \im f(1) \leq f(r) \im f(1)\in F_d$.
 Hence, for $i_0=1$ we have $\neg f(r) \vee f(i_0) \in F_d$. Besides, $\neg g(i_0 + q \frac{n}{d}-1) \vee    g(r + ( q-1) \frac{n}{d})) = \neg g(q \frac{n}{d}) \vee    g(r + ( q-1) \frac{n}{d})) \geq  \neg g(q \frac{n}{d}) \vee    g(1 + ( q-1) \frac{n}{d})) \in F_d$. So   $\xi_{d,q}(f \deb g) \in F_d$.

If $r > \frac{n}{d}$ then $\frac{n}{d} +1 \leq r \leq (d-q+1) \frac{n}{d}$ from we have $2 \leq q_{d,r} \leq d-q+1$ and $1 \leq q_{d,r} +q -1 \leq d$. Since $\neg f(r) \vee f(( q_{d,r} -1) \frac{n}{d} +1) \geq \neg f(q_{d,r} \frac{n}{d}) \vee f(( q_{d,r} -1) \frac{n}{d} +1) \in F_d$, being $i_0= ( q_{d,r} -1) \frac{n}{d} +1$ we obtain $\neg f(r) \vee f(i_0) \in F_d$. It is easy to see that $1 < i_0 \leq n- q  \frac{n}{d} +1$. Moreover, $\neg g(i_0 + q \frac{n}{d}-1) \vee    g(r + ( q-1) \frac{n}{d})) = \neg g((q_{d,r} +q -1 ) \frac{n}{d}) \vee    g(r + ( q-1) \frac{n}{d}) \geq \neg g((q_{d,r} +q -1 ) \frac{n}{d}) \vee    g((q_{d,r} +q -2 ) \frac{n}{d} +1 ) \in F_d$. Therefore $\xi_{d,q}(f \deb g) \in F_d$.\hfill$\mbox{\boldmath{$\blacksquare$}}$
\end{itemize}

\rr

In what follows let $n\geq 1$ be an integer and let $A$ be an $(n+1)-$valued Wajsberg algebra. 
For each $d \in Div(n)$, let $\chi_d(A)= \{P \in \chi (A): A/P \approx L_{d+1}\}$. 

\rr

\begin{sob}\label{O2} {\rm If $\chi_d (A) \not= \emptyset$ then there exists an element $x_d \in A$ which verifies \linebreak $g(x_d)=g(x_d \vee \neg x_d^{d-1})= \frac{d-1}{d}$, for all homomorphism $g: A \fun L_{n+1}$ such that $g(A) \approx L_{d+1}$. 

In fact, suppose that there exists $P \in \chi_d (A)$. Let $Q_d= \bigcap \{P \in \chi  (A):  A/P \approx L_{d+1}\}$. 
Let $\equiv_{Q_d}$ and $\equiv_{P}$ be congruences determined by $Q_d$ and $P$, respectively. 
It is clear that $Q_d \subseteq P$ for every $P \in  \chi_d (A)$. The following statements hold:

\begin{itemize}

\item[{\rm (2.1)}] (Theorem 6.15, [2]) The application $\alpha: (A/Q_d)/(P/Q_d) \fun A/P$ defined by the stipulation $\alpha ( [ [ a ]_{Q_d} ]_{P/Q_d})= [a]_P$  is an isomorphism  from $(A/Q_d)/(P/Q_d)$ into $A/P$, where $\equiv_{P/Q_d} = \{ ([a]_{Q_d}, [b]_{Q_d}): (a,b) \in \equiv_{P}\}$.

\item[{\rm (2.2)}] $A/Q_d$ is isomorphic to $L_{d+1}^{\chi_d(A)}$.

For each $P \in \chi_d(A)$, let $\xi_P$ be the isomorphism from  $A/P$ into $L_d+1$; besides from (2.1), there exists an isomorphism $\alpha_P : (A/Q_d)/(P/Q_d) \fun A/P$. Thus, we can consider the application $\xi_P \circ \alpha_P \circ q$ from $A/Q_d$ into $L_{d+1}$, for each $P \in \chi_d(A)$, where $q$ is the natural map from $A/Q_d$ onto $(A/Q_d)/(P/Q_d)$. Let $\psi: A/Q_d \fun  L_{d+1}^{\chi_d(A)}$ defined in the following way: 
for each $[a]_{Q_d} \in A/Q_d$, $\psi ([a]_{Q_d} )$ is the function from $\chi_d(A)$ into  $L_{d+1}$ defined by the stipulation $\psi ([a]_{Q_d} ) (P)= (\xi_P \circ \alpha_P \circ q) ([a]_{Q_d})$.

It is clear that $\psi$ is well defined and it is an homomorphism. Besides $\psi$ is bijective. In fact, 
if $[a]_{Q_d} \not= [b]_{Q_d}$ then $a \not\equiv_{Q_d} b$, i.e., there is some prime implicative filter $P \in \chi_d (A)$ such that $a \not\equiv_{P} b$. So $[a]_{Q_d} \not\equiv_{P/Q_d} [b]_{Q_d}$, and thus $q([a]_{Q_d}) \not= q([b]_{Q_d})$. Then $\psi ([a]_{Q_d}) (P) \not = \psi ([b]_{Q_d})(P)$ and $\psi ([a]_{Q_d})\not = \psi ([b]_{Q_d})$. Therefore $\psi$ is injective.

Let $h : \chi_d(A) \fun L_{d+1}$. For all $P \in \chi_d(A)$ we have $\alpha_P^{-1} ( \xi_P^{-1} (h(P))) \in (A/Q_d)/(P/Q_d)$, thus there is $a \in A$ such that $q([a]_{Q_d})=\alpha_P^{-1} ( \xi_P^{-1} (h(P)))$. It is easy to see that $\psi([a]_{Q_d}) (P)= h(P)$ for each $P \in \chi_d (A)$, i.e., $\psi([a]_{Q_d}))=h$. Therefore $\psi$ is surjective.

\item[{\rm (2.3)}] There exists an element $x_d \in A$ which verifies $h(x_d)= \frac{d-1}{d}$ for all homomorphism \linebreak $h: A \fun L_{n+1}$ such that $h(A) \approx L_{d+1}$. It suffices to take  $x_d \in A$ such that  $\psi([x_d]_{Q_d})$ is the constant function equal to $\frac{d-1}{d}$ and to note that $[x_d]_{P}=[x_d \vee \neg x_d^{d-1}]_{P}= \frac{d-1}{d}$, for every $P \in \chi_d (A)$.
\end{itemize}
 }
\end{sob}

\begin{sob}\label{O3} {\rm For each $x \in A$, any $d \in Div(n)$ and all homomorphism $g: A \fun L_{n+1}$ is $g(x \vee \neg x^{d-1}) \geq \frac{d-1}{d}$.

Indeed, suppose $g(x)= \frac{a}{n}$, for some integer $a$, $0 \leq a \leq n$. 
 If $g(x) \geq \frac{d-1}{d}$ then $g(x \vee \neg x^{d-1}) \geq \frac{d-1}{d}$. 
If $g(x) < \frac{d-1}{d}$ then $\frac{d-1}{d} < (d-1) \frac{n-a}{n}$. 
 On the other hand, $\neg (g(x))^{d-1}= (d-1) \neg g(x)= (d-1)\frac{n-a}{n}$. Hence, $\neg (g(x))^{d-1}= 1$ or $\neg (g(x))^{d-1}=  (d-1)\frac{n-a}{n} > \frac{d-1}{d}$. Thus $g(x \vee \neg x^{d-1}) \geq \frac{d-1}{d}$.}
\end{sob}

\begin{steo}\label{Teo2}
Let $\langle A, \im , \neg , 1 \rangle$ be an $(n+1)-$valued Wajsberg algebra. For each $d \in Div(n)$ let $h_A(d)= P_d \cap B(A)$, where $P_d= \bigcap \{P \in \chi  (A):  A/P \subseteq L_{d+1}\}$. Then $\varphi: A \fun M(B(A),h_A)$ is a $W-$isomorphism, being $\varphi (x)(i)=\sigma_i(x)$ for all $x \in A$ and every integer $1 \leq i \leq n$.
\end{steo}

\noindent {\bf Proof.} If $n=1$ it is trivial. Suppose $n>1$.

(i) $\varphi$ is well defined. It is clear that $\xi_{d,q}(\varphi (x)) = \neg \sigma_{q \frac{n}{d}} (x)  \vee \sigma_{(q-1) \frac{n}{d}+1} (x)  \in B(A)$ for every $d \in Div(n)$ and each integer $1 \leq q \leq d$. We shall prove that $\xi_{d,q}(\varphi (x)) \in P_d$. 

Suppose that there are $d \in Div(n)$ and an integer $1 \leq q \leq d$ such that  $\xi_{d,q}(\varphi (x)) \notin P_d$. Then there exists $P \in \chi (A)$ which verifies $A/P \approx L_{d+1}$ and $\xi_{d,q}(\varphi (x)) \notin P$. Let $[x]_P= \frac{a}{d}$, for some integer $0 < a < d$ (because if $[x]_P  \in \{0,1\}$ it is clear that $\xi_{d,q}(\varphi (x)) \in P$). 
Since $[ \xi_{d,q}(\varphi (x))]_P = \neg \sigma_{q \frac{n}{d}} [x]_P  \vee \sigma_{(q-1) \frac{n}{d}+1} [x]_P  \in \{0,1\}$, we have $[ \xi_{d,q}(\varphi (x))]_P =0$, i.e., $\sigma_{(q-1) \frac{n}{d}+1} [x]_P =0$ and $\sigma_{q \frac{n}{d}} [x]_P=1$. 
Besides, $\sigma_{q \frac{n}{d}} [x]_P=1$ if and only if $\sigma_{q \frac{n}{d}} (\frac{a \frac{n}{d}}{n})=1$ if and only if  $(q+a) \frac{n}{d} > n$ if and only if $d < a+q$ (13). 
On the other hand, $\sigma_{(q-1) \frac{n}{d}+1} [x]_P =0$ if and only if $(q-1) \frac{n}{d} +1 + a  \frac{n}{d} \leq n$ if and only if $a+q \leq d -  \frac{d}{n} +1$. Since $a + q \geq 2$ it results $1 \leq a +q -1 \leq d \frac{n-1}{n}< d$ (14). 
From (13) and (14) we obtain $d \leq q+a-1 < d$ which is a contradiction.

(ii)  $\varphi$ is a $W-$homomorphism. It is immediate because the implication is determined by the order.

(iii)  $\varphi$ is injective. It is immediate from the Moisil's determination principle. 

(iv)  $\varphi$ is surjective.  Let $f \in M\left(B(A), h_A \right)$, i.e., $f(1) \leq f(2)\leq \cdots  \leq f(n)$ and $\xi_{d,q}(f) = f(q \frac{n}{d}) \im f((q-1) \frac{n}{d}+1)  \in P_d \cap B(A)$ for all $d \in D = Div(n)$ and every integer $q$, $1 \leq q \leq d$.

We construct $z$ as follows.

For each $d \in Div(n)$, we consider the blocks corresponding to the divisor $d$. There are $d$ blocks, each one with $n/d$ elements. For each integer $i$, $1 \leq i \leq n$, the block corresponding to the divisor $d$ which contains $i$ is that determined by the integer $q_{d,i}$, the first element of the set $\{q \in \mathbb N : 1 \leq q \leq d, \, i \leq q \frac{n}{d}\}$. 
Observe that $q_{1,i}=1$ for every integer $i$, $1 \leq i \leq n$ and $q_{d,1}=1$ for all $d \in Div(n)$.

We define $y_d \in A$ as follow, where the element $x_d \in A$ is that obtained according to the Remark \ref{O2}:

$$y_d = \left\{ \begin{array}{ll}
                1 & \mbox{ if  $\chi_d (A) = \emptyset$} \\
		x_d \vee \neg x_d^{d-1} & \mbox{ in another case.}
		\end{array}
	\right.$$

For every integer $i$, $1 \leq i \leq n$, let $a_i= \bigwedge\limits_{d \in Div(n)} y_d^{(q_{d,i}-1)}$. 
Let $z \in A$, $z = \bigvee\limits_{i=1}^{n} (f(i) \wedge a_i)$. Since $a_1=1$ we can write $z = f(1) \vee \bigvee\limits_{i=2}^{n} (f(i) \wedge a_i)$. We shall prove that $\varphi (z) (j)=f(j)$ for all integer $j$, $1 \leq j \leq n$, where $\varphi (z) (j)= \sigma_j (z)$.

\rr

Suppose $f(j) \not\leq \sigma_j (z)$; then there exists an homomorphism $g: A \fun L_{n+1}$ such that $g(f(j))=1$ (15) and $g (\sigma_j (z))=0$. But $g (\sigma_j (z))= g(f(1)) \vee \bigvee\limits_{i=2}^{n} (g(f(i)) \wedge \sigma_j ( g(a_i)))=0$ if and only if (16) $g(f(1))=0$ and $g(f(i)) \wedge \sigma_j ( g(a_i))=0$ for all integer $i$, $2 \leq i \leq n$. From (15) and (16) it results $\sigma_j ( g(a_j))=0$ (17), being $\sigma_j (g(a_j))= \bigwedge\limits_{d \in D} \sigma_j \left( (g(y_d))^{(q_{d,j}-1)}\right)$.

From $q_{d,j}-1 \leq d-1 < d$ and Remark \ref{O3} we have that $(g(y_d))^{(q_{d,j}-1)}= \neg (q_{d,j}-1) \neg g(y_d) \geq \neg  (q_{d,j}-1) \frac{1}{d}= \frac{d- q_{d,j}+1}{d}$. 

Thus, $\sigma_j \left( (g(y_d))^{(q_{d,j}-1)}\right) \geq \sigma_j \left( \frac{d- q_{d,j}+1}{d} \right)= \left\{ \begin{array}{ll}
      1 & \mbox{ if  $j+ (d - q_{d,j} +1) \frac{n}{d} > n,$} \\
      0& \mbox{ in another case.}
		\end{array}
	\right.$

\rr

If $q_{d,j}=1$ then $\sigma_j \left( (g(y_d))^{(q_{d,j}-1)}\right) = 1$. Let $q_{d,j} > 1$. If $\sigma_j \left( \frac{d- q_{d,j}+1}{d} \right)= 0$ then $j \leq (q_{d,j}-1) \frac{n}{d}$ which is a contradiction because $q_{d,j}$ is the first element of the set $\{q \in \mathbb N : 1 \leq q \leq d, \, j \leq q \frac{n}{d}\}$.

Hence $\sigma_j \left( (g(y_d))^{(q_{d,j}-1)}\right) \geq \sigma_j \left( \frac{d- q_{d,j}+1}{d} \right)= 1$ for any $d \in D$, i.e., $\sigma_j ( g(a_j))=1$ which contradicts (17). Therefore $f(j) \leq \sigma_j (z)$.

\rr

Suppose now $\sigma_j (z) \not\leq f(j)$; then there exists an homomorphism $g: A \fun L_{n+1}$ such that $g (\sigma_j (z))=1$ and $g(f(j))=0$ (18). Observe that $g (\sigma_j (z))= g(f(1)) \vee \bigvee\limits_{i=2}^{n} (g(f(i)) \wedge \sigma_j ( g(a_i)))=1$ if and only if there is an integer $i_0$, $2 \leq i_0 \leq n$, which verifies $g(f(i_0)) \wedge \sigma_j ( g(a_{i_0}))=1$, if and only if $g(f(i_0))=1$ (19) and $\sigma_j ( g(a_{i_0}))=1$. Then, 
$\sigma_j (g(a_{i_0}))= \linebreak  \bigwedge\limits_{d \in D} \sigma_j \left( (g(y_d))^{(q_{d,i_0}-1)}\right)=1$ if and only if $\sigma_j \left( (g(y_d))^{(q_{d,i_0}-1)}\right)=1$ for all $d \in D$. In par\-ti\-cu\-lar, $g(A)$ is isomorphic to $L_{d_0 +1}$ for some $d_0 \in Div(n)$, hence $\sigma_j \left( (g(y_{d_0}))^{(q_{d_0,i_0}-1)}\right)=1$. If $d_0=1$ then $f(n) \im f(1)=0 \in P_1 \cap B(A)$ which is a contradiction. Thus, $d_0 >1$. By Remark \ref{O3} we have that  $(g(y_{d_0}))^{(q_{d_0,i_0}-1)}= \neg  (q_{d_0,i_0}-1) \frac{1}{d_0}= \frac{d- q_{d_0,i_0}+1}{d_0}$. Then 
$1= \sigma_j \left(    (g(y_{d_0}))^{(q_{d_0,i_0}-1)}\right)=  \sigma_j \left(   \frac{(d- q_{d_0,i_0}+1) \frac{n}{d_0}}{n}    \right)$ if and only if $j+(d- q_{d_0,i_0}+1) \frac{n}{d_0} >n$, if and only if $(q_{d_0,i_0} - 1) \frac{n}{d_0} +1 \leq j$ (20).

As $(q_{d_0,i_0} - 1) \frac{n}{d_0} +1 \leq i_0 \leq q_{d_0,i_0} \frac{n}{d_0}$ from (19) it follows $g \left( f \left( q_{d_0,i_0} \frac{n}{d_0} \right) \right)=1$. Besides $g \left( f \left( q_{d_0,i_0} \frac{n}{d_0} \right) \right) \im g \left( \left( (q_{d_0,i_0} - 1) \frac{n}{d_0} +1 \right) \right)=1$ because $\xi_{d_0,q}(f) \in P_{d_0} \cap B(A)$, then by (20) we obtain $g(f(j))=1$ which contradicts (18). Therefore $\sigma_j (z) \leq f(j)$.\hfill$\mbox{\boldmath{$\blacksquare$}}$

\rr

\begin{sdefi} {\rm (a)} A pair $\langle B,h \rangle \in B^{n+1}$ if $B$ is a Boolean algebra and $h$ is a function from the lattice of divisors of $n$ into the lattice of filters of $B$ such that $h(n)=\{1\}$ and $h(gcd \{d,r\})= h(d) \vee h(r)$, for every $d,r \in Div (n)$ {\rm (}$gcd \{d,r\}$ is the greatest common divisor of the set $\{d,r\}${\rm )}.

{\rm (b)} Objects $\langle B_1,h_1 \rangle$ and $\langle B_2,h_2 \rangle$ in $B^{n+1}$ are isomorphic if exists a boolean isomorphism $\varphi: B_1 \fun B_2$ which verifies $\varphi^{-1}(h_2(d))= h_1(d)$ for all $d \in Div (n)$.
\end{sdefi}

\begin{sob}\label{O4} {\rm Let $\langle A, \im , \neg , 1 \rangle$ be an $(n+1)-$valued Wajsberg algebra. Then $\langle B(A),h_A \rangle \in B^{n+1}$, where $h_A(d)= P_d \cap B(A)$ being $P_d= \bigcap \{P \in \chi  (A):  A/P \subseteq L_{d+1}\}$, for each $d \in Div(n)$.

In fact, let $r,t \in Div(n)$ and let $m$ be the greatest common divisor of the set $\{r,t\}$. It is clear that if $r$ divides $t$ then $P_t \subseteq P_r$, this implies $P_t \vee P_r \subseteq P_m$. For the other inclusion it is easy to verify that $A/P \subseteq L_{t+1}$ and $A/P \subseteq L_{r+1}$ implies $A/P \subseteq L_{m+1}$, for every prime implicative filter $P$ of $A$.
}
\end{sob}

\begin{steo}\label{Teo3} Let $\langle B, h  \rangle \in B^{n+1}$ and let $A=M(B,h)$. Then $\langle B, h \rangle$ and $\langle B(A), h_A \rangle$ are isomorphic objects in $B^{n+1}$.
\end{steo}

\begin{dem} Let $\langle B, h  \rangle \in B^{n+1}$ and let $A=M(B,h)$. By Corollary \ref{C1} we know that $\langle M(B,h), \linebreak \deb, \neg, \mathbb I \rangle$  is an $(n+1)-$valued Wajsberg algebra where $\mathbb I$ is the constant function equal to $1$, $(\neg f)(k)= \neg f(n+1-k)$ and $(f \deb g)(k)= \bigwedge\limits_{i=1}^{n-k+1} (f(i) \im g(i+k-1))$, for all $f,g \in M(B,h)$ and every integer $1 \leq k \leq n$. It is easy to see that $B(A)$ is the subalgebra that consist of all constant functions. If $h_A$ is the function from the lattice of divisors of $n$ into the lattice of filters of $B$ defined by $h_A (d)=P_d \cap B(A)$, being $P_d= \bigcap \{P \in \chi  (A):  A/P \subseteq L_{d+1}\}$, then $\langle B(A), h_A \rangle \in B^{n+1}$ (because Remark \ref{O4}).

Let $\mu: B \fun B(A)$ such that $\mu (a)$ is the constant function from $\{1,2, \ldots , n\}$ into $B$ that takes the value $a$, for each $a \in B$. It is clear that $\mu$ is a boolean isomorphism from $B$ onto $B(A)$. To complete the proof we should prove that $\mu^{-1}(P_d \cap B(A))=h(d)$, for each $d \in Div(n)$.

There exists isomorphisms $\psi_1$ from $\chi(B)$ onto $\chi(B(A))$ and $\psi_2$ from $\chi(B(A))$ onto $\chi (A)$ defined in the following way. For each ultrafilter $P$ of $B$, $\psi_1 (P)=\mu(P)=\{ \mu (x): x \in P \}$ and for each ultrafilter $Q$ of $B(A)$, $\psi_1^{-1} (Q)=\{a \in B: \mu (a) \in P \}$ y $\psi_2 (Q)=\{x \in A: x^n \in Q\}$. Besides, $\psi_2^{-1} (P)= P \cap B(A)$ for each prime implicative filter $P$ of $A$. 

Therefore, given a prime implicative filter $P$ of $A=M(B,h)$, we consider $Q=\psi_1^{-1}\psi_2^{-1}(P)$; $Q$ is an ultrafilter of $B$. The following statements are true:

\begin{itemize}
\item[{\bf Fact 1:}] $B/Q$ is a simple Boolean algebra and $\widetilde{h} (d) = \{ [x]_Q : x \in h(d) \}$ is a filter of $B/Q$ for every $d\in Div(n)$.

\item[{\bf Fact 2:}] $M(B/Q, \widetilde{h} (d))$ is an $(n+1)-$valued Wajsberg algebra. 
The function $\alpha$ from \linebreak $M(B/Q, \widetilde{h} (d))$ into $L_{n+1}$ defined by $\alpha (f)=\frac{k}{n}$, where $k$ is the number of in\-de\-xes in which the function $f$ takes the value $1$, is an injective homomorphism. It is clear that $\alpha$ is an isomorphism if $M(B/Q, \widetilde{h} (d))= (B/Q)^{[n]}$.

\item[{\bf Fact 3:}] $A/P$ and $M(B/Q, \widetilde{h} (d))$ are isomorphic algebras.

Indeed, let $\eta : A/P \fun B/Q$ defined for all $f \in A/P$ as $\eta ( [f]_P) (i) = [ f(i)]_Q$, for each  $i \in \{1, \ldots , n\}$. 

Let $f,g \in A=M(B,h)$ such that $f$ and $g$ are $P-$congruent. Then there exists $a \in P$ such that $f \wedge a = g \wedge a$, i.e., $f(i) \wedge a(i) = g(i) \wedge a(i)$, for each $i \in \{1, \ldots , n\}$. Since $a^n \in P \cap B(A)= \psi_2^{-1}(P)$ and $a^{n}(i)=a(1)$ for each $i$, we have that $a(1)= \mu^{-1}(a^n) \in \psi_1^{-1}\psi_2^{-1}(P)=Q$; besides $f(i) \wedge a(1) = g(i) \wedge a(1)$, thus $f(i)$ is $Q-$congruent with $g(i)$ for every $i \in \{1, \ldots , n\}$. Therefore $\eta$ is well defined.

Let $f,g \in A=M(B,h)$ such that $\eta ( [f]_P) = \eta ( [g]_P)$, then $[ f(i)]_Q = [ g(i)]_Q$ for all $i \in \{1, \ldots , n\}$. Thus there are elements $a_1, a_2, \ldots , a_n \in Q$ which verify $f(i) \wedge a_i = g(i) \wedge a_i$, for every $i \in \{1, \ldots , n\}$. Let $a = a_1 \wedge a_2 \wedge \ldots \wedge a_n \in Q$; it is clear that $\mu(a) \in P$ and $f \wedge \mu(a) = g \wedge \mu (a)$, i.e., $[f]_P = [g]_P$, hence $\eta$ is injective. It is easy to see that $\eta$ is an isomorphism.
\end{itemize}

We shall prove now $\mu^{-1}(P_d \cap B(A)) = h(d)$.

i) $\mu^{-1}(P_d \cap B(A)) \subseteq h(d)$.

Let $a \in \mu^{-1}(P_d \cap B(A))$, then $\mu (a)\in P_d \cap B(A)$, i.e., $\mu (a)\in P\cap B(A)$ for every prime implicative filter $P$ of $A$ such that $A/P \subseteq L_{d+1}$. 

If $a \notin h(d)$ there is a prime filter $Q$ of $B$ such that $h(d) \subseteq Q$ and $a \notin Q$. From Fact 3 we have that $A/P_0 \approx M(B/Q, \widetilde{h} (d))$, where  $P_0=\psi_2 \psi_1 (Q)$. Since $a \notin h(d)$ and $[a]_Q=0$ it results $\widetilde{h} (d)= \{1\}$. Hence each element in $M(B/Q, \widetilde{h} (d))$ has the {\em blocks} corresponding to the divisor $d$ only composed by $0$ or only composed by $1$; each {\em block} corresponding to the divisor $d$ have $n/d$ elements, thus the quantity of $1$ that has an element of $M(B/Q, \widetilde{h} (d))$  is a multiple of $n/d$. Hence from Fact 2, if $f \in M(B/Q, \widetilde{h} (d))$ has $k$ components equal to $1$, then $\alpha (f)=k/d$, i.e., $A/P_0 \approx M(B/Q, \widetilde{h} (d)) \subseteq L_{d+1}$. Therefore $\mu (a) \in \psi_2^{-1}(P_0)$, i.e., $a \in Q$ which is a contradiction.

ii) $h(d) \subseteq \mu^{-1}(P_d \cap B(A))$.

Let $a \in h(d)$. If $ a \notin \mu^{-1}(P_d \cap B(A))$ then there exists a prime implicative  filter $P$ of $A$ such that $A/P \subseteq L_{d+1}$ and $\mu(a) \notin P \cap B(A)$. Thus $a \notin \psi_1^{-1} \psi_2^{-1}(P)=Q_0$ ($Q_0$ is maximal among all the filters of $B$ which not containing the element $a$). It is clear that $a \not= 1$, then it exists at least a positive integer $r \in Div(n)$ such that $a \notin h(r)$. Let $\{r_1, \ldots , r_s\}= \{r \in Div(n): a \notin h(r)\}$ and let $m$ be greatest common divisor of the set $\{r_1, \ldots , r_s\}$. Since $a \notin h(r_j)$ we have $h(r_j) \subseteq Q_0$ for every $j$, $1 \leq j \leq s$. Hence $h(m)= h(r_1) \vee h(r_2) \vee \ldots \vee h(r_s) \subseteq Q_0$, from it follows $a \notin h(m)$. 
On the other hand, $M(B/Q_0, \widetilde{h} (d))$ contains a copy of $L_{m+1}$. In fact, let $S$ be the set of all $n-$tuples composed by $m$ {\em blocks}, each one of them with $\frac{n}{m}$ elements, all them equal to $0$ or all them equal to $1$. It is clear that $S \subseteq M(B/Q_0, \widetilde{h} (d))$; indeed, denote by $\widetilde{\frac{m}{n}}$ an element  of $S$. Let $r \in Div (n)$ and let $q$ be an integer between $1$ and $r$. If $a \in h(r)$ then $[a]_{Q_0}=0 \in  \widetilde{h} (r)$, therefore $ \widetilde{h} (r)= \{0,1\}$. If $a \notin h(r)$ then $m$ divides $r$; thus the {\em blocks} by $r$ are contained into the {\em blocks} by $m$; hence  $\xi_{r,q}(\widetilde{\frac{m}{n}})= 1 \in \widetilde{h} (r)$. Thus, from Fact 2 we can affirm that $L_{m+1} \subseteq M(B/Q_0, \widetilde{h} (d)) \approx A/P \subseteq L_{d+1}$, i.e., $m$ divides $d$, therefore $h(d) \subseteq h(m)$ and $a \in h(m)$ which is a contradiction.
\end{dem}

\rr

Let $\mbox{\boldmath{${\cal W}^{n+1}$}}$ be the category of $(n+1)-$valued $W-$algebras and $W-$homomorphisms. Let $\mbox{\boldmath{${\cal B}^{n+1}$}}$ be the category whose objects are pairs in $B^{n+1}$ and whose morphisms are defined in the following way: if $O_1= \langle B_1,h_1 \rangle$ and $O_2= \langle B_2,h_2 \rangle$ are objects in this category, $\theta$ is a morphism from $O_1$ into $O_2$ if it is a boolean homomorphism from $B_1$ into $B_2$ which verifies $h_1(d) \subseteq \theta^{-1}(h_2(d))$ for any $d \in Div(n)$. 

It is easy to see that $\theta$ is an isomorphism from $O_1$ onto $O_2$ if it is a boolean  isomorphism from $B_1$ onto $B_2$ which verifies $h_1(d) = \theta^{-1}(h_2(d))$ for each $d \in Div(n)$. 

Let $\mbox{\boldmath{$B$}}$ be defined from $\mbox{\boldmath{${\cal W}^{n+1}$}}$ to the category $\mbox{\boldmath{${\cal B}^{n+1}$}}$ as follows:

(i) For each object ${\cal A} = \langle A, \im, \neg, 1 \rangle$ in the category  $\mbox{\boldmath{${\cal W}^{n+1}$}}$, $\mbox{\boldmath{$B$}}({\cal A})= \langle B(A), h_A \rangle$, where $B(A)$ is the set of boolean elements of $A$ and for all $d$ divisor of $n$, $h_A (d)=P_d \cap B(A)$, being $P_d= \bigcap \{P \in \chi  (A):  A/P \subseteq L_{d+1}\}$.

(ii) If ${\cal A}_1$ and ${\cal A}_2$ are objects in the category $\mbox{\boldmath{${\cal W}^{n+1}$}}$ and $g: {\cal A}_1 \longrightarrow {\cal A}_2$ is a $\mbox{\boldmath{${\cal W}^{n+1}$}}-$mor\-phism, $\mbox{\boldmath{$B$}}(g):  \langle B(A_1), h_{A_1} \rangle \longrightarrow \langle B(A_2), h_{A_2} \rangle$  is defined by $\mbox{\boldmath{$B$}}(g)= g/_{B(A_1)}$.

It is immediate that $\mbox{\boldmath{$B$}}(g)$ is a boolean homomorphism. Moreover, $\mbox{\boldmath{$B$}}(g)$ is a $\mbox{${\cal B}^{n+1}$}-$mor\-phism. Indeed, let $a \in h_{A_1} (d)$. If $a \notin \mbox{\boldmath{$B$}}(g)^{-1} (h_{A_2} (d))$ then $g(a) \notin h_{A_2} (d)$, hence there exists a prime implicative filter $P$ of $A_2$ such that $A_2/P \subseteq L_{d+1}$ and $g(a) \notin P$. Thus $a \notin g^{-1}(P) \cap B(A_1)$. The function $v: A_1/g^{-1}(P) \fun A_2/P$ defined by $v([x]_{g^{-1}(P)})= [g(x)]_P$ is an embedding from $A_1/g^{-1}(P)$ into $A_2/P \subseteq L_{d+1}$, i.e., $A_1/g^{-1}(P) \subseteq L_{d+1}$ then $a \notin h_{A_1} (d)$ which is a contradiction. 

It is easy to verify that $\mbox{\boldmath{$B$}}$ is a functor.

\rr

Let $\mbox{\boldmath{$M$}}$ be defined from $\mbox{\boldmath{${\cal B}^{n+1}$}}$ to  $\mbox{\boldmath{${\cal W}^{n+1}$}}$ as follows:

(i) For each object $\langle B, h \rangle$ in the category $\mbox{\boldmath{${\cal B}^{n+1}$}}$, let $\mbox{\boldmath{$M$}} (\langle B, h \rangle)= \langle M(B,h), \deb, \neg, \mathbb I \rangle$.

(ii) If $\langle B_1, h_1 \rangle$ and $\langle B_2, h_2 \rangle$ are objects in the  category  $\mbox{\boldmath{${\cal B}^{n+1}$}}$ and $g$ is a $\mbox{\boldmath{${\cal B}^{n+1}$}}-$morphism from $\langle B_1, h_1  \rangle$ into $\langle B_2, h_2 \rangle$ let $\mbox{\boldmath{$M$}}(g): M( B_1, h_1)  \longrightarrow M( B_2, h_2)$ where $\mbox{\boldmath{$M$}}(g) (f)= g \circ f$, for any $f \in M( B_1, h_1)$. 

It is clear that $\mbox{\boldmath{$M$}}(g)$ is well defined because, if $f \in M( B_1, h_1)$ then for each $d \in Div(n)$ and all integer $q$, $1 \leq q \leq d$ we have $\xi_{d,q}(f) \in h_1(d)$; hence $\xi_{d,q}(g \circ f)= g(\xi_{d,q}(f)) \in g(h_1(d)) \subseteq g g^{-1} (h_2(d) \subseteq h_2(d)$. Therefore $g \circ f \in M( B_2, h_2)$. Besides $\mbox{\boldmath{$M$}}(g)$ is a $W-$homomorphism.

It is verifies that $\mbox{\boldmath{$M$}}$ is a functor.

From Theorems \ref{Teo2} and \ref{Teo3} it is easy to see that the functors $\mbox{\boldmath{$B$}}$ and $\mbox{\boldmath{$M$}}$ define a natural equivalence between the categories $\mbox{\boldmath{${\cal W}^{n+1}$}}$ and $\mbox{\boldmath{${\cal B}^{n+1}$}}$.

\rr

\rr

In [5] R. Cignoli defines the $(n+1)-$valued Boolean spaces and establishes the result which appears below as Theorem \ref{TEORC}.

\begin{sdefi} {\rm (i)} An $(n+1)-$valued Boolean space is a pair $\langle X,h \rangle$ where $X$ is a Boolean space and $h$ is a meet-homomorphism from the lattice of positive divisors of $n$ into the lattice of closed subsets $X$, such that $h(n)=X$.

{\rm (ii)} A morphism from $\langle X,h \rangle$ into $\langle Y,g \rangle$ is a continuous function  $\varphi: X \fun Y$ which verifies $\varphi^{-1}(g(d))= h(d)$ for each divisor $d$ of $n$. 
\end{sdefi}

\begin{steo}\label{TEORC} {\rm [5]} For each $(n+1)-$valued $MV-$algebra $A$, there is a unique (up to  isomorphisms) $(n+1)-$valued Boolean space $\langle X,h \rangle$ such that $A$ is isomorphic to  the $(n+1)-$valued $MV-$algebra formed by the continuous functions $f: X \fun L_{n+1}$ such that $f(h(d)) \subseteq L_{d+1}$, for each divisor $d$ of $n$, with the $MV-$operations defined pointwise. The space $X$ is homeomorphic to the Stone space of the maximal Boolean subalgebra of $A$.
\end{steo}

\begin{sob} \label{O5}{\rm Theorems \ref{Teo2} and \ref{TEORC} are  equivalent in the sense that each one of them can be deduced from the other one. The unicity in the Theorem \ref{TEORC} follows from Theorem \ref{Teo3}.

In fact, let $A$ be an $(n+1)-$valued $MV-$algebra. If Theorem \ref{Teo2} is true we know that $A$ es isomorphic to  $\M=M(B(A),h^{\prime})$, where $h^{\prime}$ is a function from the lattice of divisors of $n$ into the lattice of filters of $B(A)$ such that $h(n)=\{1\}$ and $h(gcd \{d,r\})= h(d) \vee h(r)$, for every $d,r \in Div (n)$, being $gcd \{d,r\}$ the greatest common divisor of the set $\{d,r\}$.
 
$\chi(B(A))$ is a topological space where the clopen sets $s(a)=\{ F \in \chi(B(A)): a \in F\}$, for each $a \in B(A)$, forms a basis for this topology.

It is know that there is an order-antisomorphism $\delta$ from the lattice of filters of $B(A)$ onto the lattice of  closed sets of $\chi(B(A))$ defined by $\delta (F)= \cap \{s(x): x \in F\}$, for each filter $F$ of $B(A)$. 
Let us consider the application $h=\delta \circ h^{\prime}$ from the lattice of divisors of $n$ into the lattice of closed subset of $\chi(B(A))$. It is clear that  $h$ is a meet-homomorphism and $h(n)=\chi(B(A))$. Thus, $\langle \chi(B(A)), h \rangle$ is an $(n+1)-$valued boolean space.
We denote by $\C$ the set of all continuous functions $f$ from $\chi(B(A))$ into $L_{n+1}$ which verify $f(h(j)) \subseteq L_{j+1}$ for each $j \in Div(n)$. Let $\Psi: \M \fun \C$ be the function defined as follows. For every $g \in \M$ and each integer $1 \leq i \leq n$,  let $S_i=s(g(i))= \{U \in \chi(B(A)): g(i) \in U \}$; then $\Psi(g)=\widetilde{g}$, where  $\widetilde{g}: \chi(B(A)) \fun L_{n+1}$ is defined by:

$$\widetilde{g} (U) = \left\{ \begin{array}{ll}
                1 & \mbox{ if  $U \in S_1$,} \\
		 \frac{k}{n} & \mbox{ if $U \in S_{n-k+1} \setminus S_{n-k}$, for $1 \leq k \leq n-1$,} \\
		0 & \mbox{ if $U \notin S_n$.}
		\end{array}
	\right.$$

It is easy to see that $\widetilde{g}$ is continuous because $L_{n+1}$ has the discrete topology, $\widetilde{g}^{-1}(\{1\})=s(g(1))$, $\widetilde{g}^{-1}(\{0\})=\chi(B(A)) \setminus s(g(n))$ and for each $1 \leq k \leq n-1$, $\widetilde{g}^{-1}(\{\frac{k}{n}\})=s(g(n-k+1)) \cap (\chi(B(A)) \setminus s(g(n-k)))$, all open sets in $\chi(B(A))$.

Let $j \in Div(n)$; for $j=n$ we have $\widetilde{g} (h(n))= \widetilde{g} (\chi(B(A)) \subseteq L_{n+1}$. For $j \not= n$, given $U \in h(j)=\delta (F_j)=\cap\{s(a): a \in F_j\}$ we have that  $U \in s(a)$ for every $a \in F_j$, i.e., $F_j \subseteq U$. If $g(1) \in U$ then  $\widetilde{g}(U)=1 \in L_j+1$. If $g(1) \notin U$, let $r$ be the greatest positive integer which verifies $g(i)\notin U$, for $1 \leq i \leq n$. If $r=n$ then $\widetilde{g}(U)=0 \in L_j+1$. If $r < n$ then $g(r+1) \in U$. Let $q$ be the first element of the set $\{t \in \mathbb N : 1 \leq t \leq j, r \leq t \frac{n}{j}\}$. Since $g \in \M$ we have $g(q \frac{n}{j}) \im g((q-1) \frac{n}{j}+1) \in F_j \subseteq U$, hence $\neg g(q \frac{n}{j}) \in U$ o $g((q-1) \frac{n}{j}+1) \in U$. If $\neg g(q \frac{n}{j}) \in U$ then $U \notin s(g(q \frac{n}{j}))$ from we obtain $r=q \frac{n}{j}$; in this case $U \in S_{r+1}\setminus S_r$ and  $\widetilde{g}(U)=1-\frac{q}{j} \in L_{j+1}$. If $g((q-1) \frac{n}{j}+1) \in U$ then $r \leq (q-1) \frac{n}{j}$; thus $q-1 \in \{t \in \mathbb N : 1 \leq t \leq j, r \leq t \frac{n}{j}\}$ which is a contradiction. Therefore $\widetilde{g} \in \C$ and $\Psi$ is well defined.

Let $f,g \in \M$ such that $\Psi (f) = \Psi (g)$. Thus $\widetilde{f} (U)= \widetilde{g} (U)$ for each $U \in \chi(B(A))$. If there exists an integer $j \in \{1, \ldots, n\}$ for the one which $f(j) \not= g(j)$ then there is $U \in \chi(B(A))$ such that $f(j) \in U$ and $g(j) \notin U$ (or $f(j) \notin U$ and $g(j) \in U$, in which case the proof is similar).  
Let $j_0$ be the greatest positive integer which verifies $g(j) \notin U$. If $j_0 = n$ then  $\widetilde{g} (U)=0$ and $\widetilde{f} (U)\not= 0$ because $U \in s(f(j)) \subseteq s(f(j_0))$. 
If $j_0 < n$ then $g(j_0) \notin U$ and $g(j_0+1) \in U$, in this case $\widetilde{g} (U)= \frac{n-j_0}{n}$ and $\widetilde{f} (U) \not= \frac{n-j_0}{n}$ because $U \in s(f(j_0))$. Therefore $\Psi$ is injective.

Let $f \in \C$, i.e., $f$ is a continuous function from $\chi(B(A))$ into $L_{n+1}$ which verifies $f(h(d)) \subseteq L_{d+1}$ for each $d \in Div(n)$. Let $A_0=\emptyset$, and for each integer $j$, $1 \leq j \leq n$ let $A_j=A_{j-1} \cup \{ \frac{n-j+1}{n}\}$. Each $A_j$ is a clopen subset of $L_{n+1}$, then $f^{-1}(A_j)$ is a clopen subset of $\chi(B(A))$, therefore   there are $a_1, a_2, \ldots, a_n \in B(A)$ such that $s(a_j)= f^{-1}(A_j)= \{U \in \chi(B(A)): a_j \in U\}$. The conditions (i) and (ii) hold.

(i) $a_i \leq a_j$ if $i \leq j$.

It is enough to prove $s(a_i) \subseteq s(a_j)$ for every $1 \leq i \leq j \leq n$. Let $U \in s(a_i)$, then $f(U) \in A_i \subseteq A_j$ from it results $U \in s(a_j)$.

(ii) For all $d \in Div(n)$ and each integer $q$, $1 \leq q \leq d$, the element $a(d,q)=a_{q \frac{n}{d}} \im a_{(q-1) \frac{n}{d}+1} \in h(d)$. 

Indeed, if $a(d,q) \notin F_d$ there exists $P_0 \in \chi(A)$ such that $A/P_0 \subseteq L_{d+1}$ and $a(d,q)= \neg a_{q \frac{n}{d}} \vee  a_{(q-1) \frac{n}{d}+1} \notin P_0$, i.e., $a_{(q-1) \frac{n}{d}+1} \notin P_0$ y $a_{q \frac{n}{d}} \in P_0$. Let $U_0= P_0 \cap B(A) \in  \chi(B(A))$. We know that $f(h(d)) \subseteq L_{d+1}$ for each $d \in Div(n)$, with  $h(d)=\delta(F_d)= \cap \{s(x): x \in F_d\}$, then it is clear that $U \in \delta(F_d)$, if and only if, $x \in F_d$ implies $x\in U$ for each $x \in B(A)$. Therefore $U_0 \in \delta (F_d)$ from we have $f(U_0) \in f(\delta (F_d)) \subseteq L_{d+1}$. Suppose $f(U_0)= \frac{k}{d}=\frac{kn/d}{n}$, with $0 \leq k \leq d$. Let $i=q \frac{n}{d}$ and $j=(q-1) \frac{n}{d}+1$, it is clear that $j \leq i$. Let $r \in \{1, 2, \ldots , n \}$ be the least integer which verifies $f(U_0) \in A_r=A_{r-1} \cup \{ \frac{n-r+1}{n}\}$. Hence $f(U_0)= \frac{n-r+1}{n}$ from we obtain $r=n-k\frac{n}{d}+1$. If $r \leq j$, $A_r \subseteq A_j$ then $a_j \in U_0$ which is a contradiction. Therefore it should be $j < r$, thus $q + k -1 <d$ (21). If $i < r$ then $f(U_0) \notin A_i$ which is impossible, hence $r \leq i$, i.e. $d \leq q+k-1$ which contradicts (21). Thus $a(d,q) \in h(d)$.

From (i) and (ii) we have that $(a_1, \ldots, a_n) \in \M$; it is easy to verify $\Psi ((a_1, \ldots, a_n))= f$, being $\Psi ((a_1, \ldots, a_n)) = \widetilde{a}$ the function from $\chi(B(A))$ into $L_{n+1}$ defined by

\rr

\noindent $\widetilde{a} (U) = \left\{ \begin{array}{ll}
                1 & \mbox{ if  $U \in s(a_1)$,} \\
		 \frac{k}{n} & \mbox{ if $U \in s(a_{n-k+1}) \setminus s(a_{n-k})$, for $1 \leq k \leq n-1$,} \\
		0 & \mbox{ if $U \notin s(a_n)$.}
		\end{array}
	\right.$

\rr

Therefore $\Psi$ is surjective.

It is immediate that Theorem \ref{Teo2} follows from Theorem \ref{TEORC} considering the function $\Psi^{-1}: \C \fun \M$.}
\end{sob} 

\

\noindent {\bf Acknowledgement.} We wish to express our gratitude to the referee for her/his valuable remarks on this paper.

\rr

\begin{center}
{\large \bf References}
\end{center}

\begin{itemize}
\item[{[1]}]  V. Boicescu, A. Filipoiu, G. Georgescu and S. Rudeanu,
{\em  {\L}ukasiewicz--Moisil \linebreak Algebras}, North--Holland, 1991.

\item[{[2]}]  S. Burris and H. P. Sankappanavar, {\em A Course in Universal Algebra}, Graduate texts in Mathematics 78, Springer Verlag, New York, 1981.

\item[{[3]}]  C.C. Chang, {\em Algebraic analysis of many valued logics}, Transactions of the American Ma\-the\-matical Society 88 (1958), 467--490.

\item[{[4]}]  C.C. Chang, {\em A new proof of the completeness of the {\L}ukasiewicz axioms}, Transactions of the American Mathematical Society 93 (1959), 74--80.

\item[{[5]}]  R. Cignoli, {\em Natural dualities for the algebras of {\LL}ukasiewicz finitely-valued logics}, The Bulletin of Symbolic Logic 2 (1996), 218.

\item[{[6]}]  R.L.O. Cignoli, I.M.L. D'Ottaviano and D. Mundici, {\em Algebraic Foundations of Many-vlued Reasoning}, Kluwer Academic Publishers, 2000. 

\item[{[7]}]  A. Di Nola and A. Lettieri, {\em One Chain Generated Varieties of MV-Algebras}, Journal of Algebra 225 (2000), 667--697.

\item[{[8]}]  J.M. Font, A.J. Rodriguez and A. Torrens, {\em Wajsberg algebras}, Stochastica 8 (1984), 5--31.

\item[{[9]}]  R.S. Grigolia, {\em Algebraic analysis of {\L}ukasiewicz-Tarski's $n-$valued logical systems}, in R. W\'ojcicki, G. Malinowski (Eds.) Selected Papers on {\L}ukasiewicz Sentential Calculi, Ossolineum, Wroc{\cal l}aw (1977), 81--92.

\item[{[10]}]  A. Iorgulescu, {\em Connections between $MV_n$ algebras and  
$n-$valued Lukasiewicz--Moisil algebras Part II}, Discrete Mathematics 202  1--3 (1999), 113--134.

\item[{[11]}]  K. Keimel and H. Werner, {\em Stone duality for varieties generated by a quasi-primal algebra}, Mem. Amer. Math. Soc. No.148 (1974), 59--85.

\item[{[12]}]  Y. Komori, {\em The separation theorem of the $\aleph_{0}-$valued {\L}ukasiewicz propositional logic}, Reports of the Faculty of Sciences, Shizuoka University  12 (1978), 1--5.

\item[{[13]}]  N. G. Mart\'{\i}nez, {\em The Priestley Duality for Wajsberg Algebras}, Studia Logica 49, No. 1 (1990), 31--46.

\item[{[14]}]  P. Niederkorn, {\em Natural Dualities for Varieties of MV-algebras, I}, Journal of Ma\-the\-matical Analysis and Applications 255 (2001), 58--73.

\item[{[15]}]  A.J. Rodr\'{\i}guez, Un estudio algebraico de los C\'alculos Proposicionales de {\L}ukasiewicz, Tesis Doctoral, Univ. de Barcelona, 1980.

\item[{[16]}]  A. J. Rodr\'{\i}guez and A. Torrens, {\em Wajsberg Algebras and Post Algebras}, Studia Logica 53 (1994), 1--19.
\end{itemize}

\

Departamento de Matem\'atica - Facultad de Ciencias Exactas y Naturales 

Universidad Nacional de La Pampa

Av. Uruguay 151 - (6300) Santa Rosa, La Pampa, Argentina. 

e-mail: {\em mblatt@exactas.unlpam.edu.ar} 

\rr

Departamento de Matem\'atica - Universidad Nacional de Buenos Aires 

Pabell\'on I - Ciudad Universitaria, Buenos Aires, Argentina.

e-mail: {\em apetrov@dm.uba.ar} 
\end{document}